\documentclass[11pt,a4paper]{article}

\usepackage{graphicx}
\usepackage{amsmath}
\usepackage{amsfonts}
\usepackage{amssymb}
\usepackage{authblk}
\usepackage[english]{babel}
\usepackage[utf8]{inputenc}
\usepackage{subfig}
\usepackage{geometry}

\title{Energy preserving moving mesh methods applied to the BBM equation}
\author{Sølve Eidnes$^*$ and Torbjørn Ringholm%
  \thanks{Department of Mathematical Sciences, Norwegian University of Science and Technology, N–7491 Trondheim, Norway\\
  Sølve Eidnes: \texttt{solve.eidnes@ntnu.no}; Torbjørn Ringholm: \texttt{torbjorn.ringholm@ntnu.no}}}
\begin{document}
\maketitle

\abstract{Energy preserving numerical methods for a certain class of PDEs are derived, applying the partition of unity method. The methods are extended to also be applicable in combination with moving mesh methods by the rezoning approach. These energy preserving moving mesh methods are then applied to the Benjamin--Bona--Mahony equation, resulting in schemes that exactly preserve an approximation to one of the Hamiltonians of the system. Numerical experiments that demonstrate the advantages of the methods are presented.}


\section{INTRODUCTION}
Numerical solutions of differential equations by standard methods will typically not inherit invariant properties from the original, continuous problem. Since the energy-preserving methods of Courant, Friedrichs and Lewy were introduced in \cite{courant28udp}, the development of conservative methods has garnered much interest and considerable research, surveyed in \cite{li95fdc} up to the early 1990s. 
In some important cases, conservation properties can be used to ensure numerical stability or existence and uniqueness of the numerical solution. In other cases, the conservation of one or more invariants can be of importance in its own right. In addition, as noted in \cite{hairer06gni}, one may expect that when properties of the continuous dynamical system are inherited by the discrete dynamical system, the numerical solution can be more accurate, especially over large time intervals.

The discrete gradient methods for ordinary differential equations (ODEs), usually attributed to Gonzalez \cite{gonzalez96tia}, are methods that preserve first integrals exactly. Since the late 1990s, a number of researchers have worked on extending this theory to create a counterpart for partial differential equations (PDEs), see e.g. \cite{furihata11dvd, celledoni12per}. Such methods, which are either called discrete variational derivative methods or discrete gradient methods for PDEs, aim at preserving some discrete approximation of a first integral which is preserved by the continuous system. Up to very recently, the schemes presented have typically been based on a finite difference approach, and exclusively on fixed, uniform grids. Two different discrete variational derivative methods on fixed, non-uniform grids were presented by Yaguchi, Matsuo and Sugihara in \cite{yaguchi10aeo, yaguchi12tdv}. In \cite{MM_Miyatake_Matsuo}, Miyatake and Matsuo introduce integral preserving methods on adaptive grids for certain classes of PDEs. Eidnes, Owren and Ringholm presented in \cite{EOR2015} a general approach to extending the theory of discrete variational derivative methods, or discrete gradient methods for PDEs, to adaptive grids, using either a finite difference approach, or the partition of unity method, which can be seen as a generalization of the finite element method.

In this paper, we present an application of the approach introduced in \cite{EOR2015} to the Benjamin--Bona--Mahony (BBM) equation, also called the regularized long wave equation in the literature. Although what we present here is a finite element method, the theory can be easily applied in a finite difference setting. Previously, there have been developed integral preserving methods for this equation \cite{cohen11gfd}, as well as adaptive moving mesh methods \cite{Lu_Huang_Qiu}, but the schemes we are to present here are, to our knowledge, the first combining these properties. In fact, in \cite{Lu_Huang_Qiu} it is noted that combining integral preservation with adaptivity is an interesting topic for further research.

\section{THE DISCRETE GRADIENT METHODS FOR PDEs}

We give a quick survey of the discrete gradient methods for PDEs, and present an approach to the spatial discretization by the partition of unity method (PUM).

\subsection{Problem statement}

Consider a PDE of the form
\begin{align}
u_t = f(\mathbf{x},u^J), \qquad \mathbf{x} \in \Omega \subseteq \mathbb{R}^d,\quad u \in \mathcal{B} \subseteq L^2, \label{eq:pde_pure}
\end{align}
where $u^J$ denotes $u$ itself and its partial derivatives of any order with respect to the spatial variables $x_1,....,x_d$, and where we assume that $\mathcal{B}$ is sufficiently regular to allow all operations used in the following.

We define a \textit{first integral} of (\ref{eq:pde_pure}) to be a functional $\mathcal{I}[u]$ satisfying
\begin{align*}
\left\langle \dfrac{\delta \mathcal{I}}{\delta u}[u], f(\mathbf{x},u^J) \right\rangle_{L^2}  = 0, \quad \forall u \in \mathcal{B},
\end{align*}
recalling that the \textit{variational derivative} $\frac{\delta \mathcal{I}}{\delta u}[u]$ is defined as the function satisfying
\begin{align*}
\left\langle \dfrac{\delta \mathcal{I}}{\delta u}[u], v \right\rangle_{L^2} = \dfrac{\mathrm{d}}{\mathrm{d}\epsilon} \bigg |_{\epsilon = 0} \mathcal{I}[u+\epsilon v] \quad \forall v \in \mathcal{B}.
\end{align*}
This means that $\mathcal{I}[u]$ is preserved over time by (\ref{eq:pde_pure}), since
\begin{align*}
\dfrac{\mathrm{d}\mathcal{I}}{\mathrm{d}t} = \left\langle \dfrac{\delta \mathcal{I}}{\delta u}[u], \dfrac{\partial u}{\partial t}\right\rangle_{L^2} = 0.
\end{align*}
Furthermore, we may observe that if there exists some operator $S(\mathbf{x},u^J)$, skew-symmetric with respect to the $L^2$ inner product, such that 
\begin{align*}
f(\mathbf{x},u^J) = S(\mathbf{x},u^J)\dfrac{\delta \mathcal{I}}{\delta u}[u],
\end{align*}
then $\mathcal{I}[u]$ is a first integral of (\ref{eq:pde_pure}), and we can state (\ref{eq:pde_pure}) on the form
\begin{align}
u_t = S(\mathbf{x}, u^J)\dfrac{\delta \mathcal{I}}{\delta u}[u]. \label{eq:pde_var}
\end{align}
The idea behind the discrete variational derivative methods is to derive a discrete version of the PDE on the form (\ref{eq:pde_var}), by obtaining a so-called discrete variational derivative 
and approximate $S(\mathbf{x},u^J)$ by a skew-symmetric matrix, see e.g. \cite{furihata11dvd}.

As proven in \cite{EOR2015}, all discrete variatonal derivative methods can be expressed as discrete gradient methods on a system of ODEs obtained by discretizing (\ref{eq:pde_var}) in space, to get a system
\begin{align}
\dfrac{\mathrm{d}\mathbf{u}}{\mathrm{d}t} = S(\mathbf{u})\nabla I(\mathbf{u}), \label{eq:ODE_form}
\end{align}
where $S(\mathbf{u})$ is a skew-symmetric matrix. The discrete gradient methods for such a system of ODEs preserve the first integral $I(\mathbf{u})$ \cite{mclachlan99giu}. These numerical methods are given by
\begin{align*}
\dfrac{\mathbf{u}^{n+1} - \mathbf{u}^n}{\Delta t} = \bar{S}(\mathbf{u}^n,\mathbf{u}^{n+1}) \overline{\nabla}I(\mathbf{u}^n,\mathbf{u}^{n+1}),
\end{align*}
where $\bar{S}(\mathbf{u}^{n},\mathbf{u}^{n+1})$ is a consistent skew-symmetric time-discrete approximation to $S(\mathbf{u})$ and $\overline{\nabla}I(\mathbf{v},\mathbf{u})$ is a discrete gradient of $I(\mathbf{u})$, defined as a function satisfying
\begin{align*}
(\overline{\nabla}I(\mathbf{v},\mathbf{u}))^T(\mathbf{u} - \mathbf{v}) &= I(\mathbf{u}) - I(\mathbf{v}),\\
\overline{\nabla}I(\mathbf{u},\mathbf{u}) &= \nabla I(\mathbf{u}).
\end{align*}
There are many possible choices of discrete gradients. For the numerical experiments in this note, we will use the Average Vector Field (AVF) discrete gradient \cite{celledoni12per}, given by
\begin{align*}
\overline{\nabla}I(\mathbf{v},\mathbf{u}) = \int \limits_0^1 \nabla I(\xi \mathbf{u} + (1 - \xi)\mathbf{v}) \mathrm{d}\xi,
\end{align*}

Note that when discretizing the system (\ref{eq:pde_var}) in space, we do so by finding a discrete approximation $\mathcal{I}_\mathbf{p}$ to the integral $\mathcal{I}$, and define an energy preserving method to be a method preserving this approximation.

\subsection{Partition of unity method on a fixed mesh}\label{sec:PUM}

The partition of unity method is a generalization of the finite element method (FEM). Stating a weak form of (\ref{eq:pde_var}), the problem consists of finding $u \in \mathcal{B}$ such that
\begin{align*}
\left\langle u_t,v \right\rangle_{L^2} = \left\langle S(\mathbf{x}, u^J)\dfrac{\delta \mathcal{I}}{\delta u}[u],v \right\rangle_{L^2} = - \left\langle \dfrac{\delta \mathcal{I}}{\delta u}[u],S(\mathbf{x}, u^J) v \right\rangle_{L^2} \quad \forall v \in \mathcal{B}.
\end{align*}
We define an approximation to $u$ by
\begin{align*}
u^h(x,t) = \sum_{i=0}^M u_i(t) \varphi_i(x),
\end{align*}
where the test functions $\varphi_i(x)$ span a finite-dimensional subspace $\mathcal{B}^h \subseteq \mathcal{B}$.
Referring to \cite{EOR2015} for details, we then obtain the Galerkin form of the problem: Find $u_i(t), i=0,\ldots,M,$ such that
\begin{align*}
\sum_{i=0}^M\dfrac{\mathrm{d} u_i}{\mathrm{d} t}\left\langle \varphi_i,\varphi_j \right\rangle_{L^2} = - \sum_{i=0}^M w_i(\mathbf{u})\left\langle \varphi_i,S(\mathbf{x}, u^{h,J}) \varphi_j \right\rangle_{L^2} \quad \forall j \in \{0,...,M\},
\end{align*}
where, with $A_{i j} = \left\langle \varphi_i, \varphi_j \right\rangle_{L^2}$,
\begin{align*}
\mathbf{w}(\mathbf{u}) = A^{-1}\nabla \mathcal{I}_{\mathbf{p}}(\mathbf{u}).
\end{align*}
We end up with an ODE for the coefficients $u_i$:
\begin{align}
\dfrac{\mathrm{d}\mathbf{u}}{\mathrm{d}t} = S_{\mathbf{p}}(\mathbf{u}) \nabla \mathcal{I}_{\mathbf{p}}(\mathbf{u}). \label{eq:PUM_ODE}
\end{align}
Here, $S_{\mathbf{p}}(\mathbf{u})  = -A^{-1} B(\mathbf{u})  A^{-1}$ is a skew-symmetric matrix, with $B(\mathbf{u})$ given by $B(\mathbf{u})_{ji} = \left\langle \varphi_i,S(\mathbf{x}, u^{h,J}) \varphi_j \right\rangle_{L^2}$, and the system is thereby of the form (\ref{eq:ODE_form}). Then, the scheme
\begin{align*}
\dfrac{\mathbf{u}^{n+1} - \mathbf{u}^n}{\Delta t} = S_{\mathbf{p}}(\mathbf{u}^{n},\mathbf{u}^{n+1}) \overline{\nabla}\mathcal{I}_{\mathbf{p}}(\mathbf{u}^{n},\mathbf{u}^{n+1}). 
\end{align*}
will preserve the approximated first integral $\mathcal{I}_{\mathbf{p}}$ in the sense that $\mathcal{I}_{\mathbf{p}}(\mathbf{u}^{n+1}) = \mathcal{I}_{\mathbf{p}}(\mathbf{u}^n)$.

\section{ADAPTIVE SCHEMES}
\label{sec:adaptive}

The primary motivation for using an adaptive mesh is usually to increase accuracy while keeping computational cost low, by improving discretization locally. Such methods are typically useful for problems with e.g. traveling wave solutions and boundary layers. The different strategies for adaptive meshes can be classified into two main groups \cite{MM_Huang_Russell}: The quasi-Lagrange approach involves coupling the evolution of the mesh with the PDE, and then solving the problems simultaneously; The rezoning approach consists of calculating the function values and mesh points in an intermittent fashion. Our method can be coupled with any adaptive mesh strategy utilizing the latter approach.

\subsection{Adaptive discrete gradient methods}

Let $\mathbf{p}^n$, $\mathbf{u}^n$, $\mathbf{p}^{n+1}$, and $\mathbf{u}^{n+1}$ denote the discretization parameters and the numerical values obtained at the current time step and next time step, respectively. Note that we now alter the notion of a preserved first integral further, to requiring that $\mathcal{I}_{\mathbf{p}^{n+1}}(\mathbf{u}^{n+1}) = \mathcal{I}_{\mathbf{p}^n}(\mathbf{u}^n)$.  The idea behind our approach is to find $\mathbf{p}^{n+1}$ based on $\mathbf{u}^n$ and $\mathbf{p}^n$, transfer $\mathbf{u}^n$ to $\mathbf{p}^{n+1}$ to obtain $\hat{\mathbf{u}}$, and then use $\hat{\mathbf{u}}$ to propagate in time to get $\mathbf{u}^{n+1}$. If the transfer operation between the meshes is preserving, i.e. if $\mathcal{I}_{\mathbf{p}^{n+1}}(\hat{\mathbf{u}}) = \mathcal{I}_{\mathbf{p}^n}(\mathbf{u}^n)$, then the next time step can be taken with the discrete gradient method for static meshes. If, however, non-preserving transfer is used, corrections are needed in order to get a numerical scheme. We introduce in \cite{EOR2015} the scheme
\begin{align}
\mathbf{u}^{n+1}\! = \! \hat{\mathbf{u}}\! - \!\dfrac{(\mathcal{I}_{\mathbf{p}^{n+1}}(\hat{\mathbf{u}}) - \mathcal{I}_{\mathbf{p}^n}(\mathbf{u}^n))\mathbf{z}}{\left\langle \overline{\nabla}\mathcal{I}_{\mathbf{p}^{n+1}}(\hat{\mathbf{u}},\mathbf{u}^{n+1}), \mathbf{z} \right\rangle}\! + \! \Delta t S_{\mathbf{p}^{n+1}}(\hat{\mathbf{u}},\mathbf{u}^{n+1}) \overline{\nabla}\mathcal{I}_{\mathbf{p}^{n+1}}(\hat{\mathbf{u}},\mathbf{u}^{n+1}),
\label{eq:DGMM}
\end{align}
where $\mathbf{z}$ is a vector which should be chosen so as to obtain a minimal correction, and such that $\langle \overline{\nabla}\mathcal{I}_{\mathbf{p}^{n+1}}(\hat{\mathbf{u}},\mathbf{u}^{n+1}), \mathbf{z} \rangle \neq 0$. In the numerical experiments to follow, we have used $\mathbf{z} = \overline{\nabla}\mathcal{I}_{\mathbf{p}^{n+1}}(\hat{\mathbf{u}},\mathbf{u}^{n+1})$.

A preserving transfer can by obtained using the method of Lagrange multipliers. Depending on whether $r$- $p$- or $h$-refinement (or a combination) is used between time steps, we expect the shape and/or number of basis functions to change. See e.g.$\!$ \cite{MM_Huang_Russell} or \cite{HPFEM_Babuska_Guo} for examples of how the basis may change through adaptivity. Denote by $\mathcal{B}^h = \mathrm{span}\{\varphi_i \}_{i=0}^M$ the trial space from the current time step and by $\hat{\mathcal{B}}^h = \mathrm{span}\{\hat{\varphi_i} \}_{i=0}^{\hat{M}}$ the trial space for the next time step, and note that in general, $M \neq \hat{M}$. We wish to transfer the approximation $u^h$ from $ \mathcal{B}^h $ to $ \hat{\mathcal{B}}^h  $, obtaining an approximation $\hat{u}^h$, while conserving the first integral, i.e. $\mathcal{I}[u^h] = \mathcal{I}[\hat{u}^h]$. This can be formulated as a constrained minimization problem:
\begin{align}
\min_{\hat{u}^h \in \tilde{\mathcal{B}^h}}||\hat{u}^h - u^h||_{L^2}^2 \quad \text{s.t.} \quad  \mathcal{I}[\hat{u}^h] = \mathcal{I}[u^h].
\label{eq:minprob}
\end{align}
Observe that
\begin{align*}
||\hat{u}^h - u^h||_{L^2}^2 &= \sum_{i=0}^{\hat{M}}\sum_{j=0}^{\hat{M}} \hat{u}_i \hat{u}_j\hat{A}_{ij} - 2\sum_{i=0}^{\hat{M}}\sum_{j=0}^M \hat{u}_i u_j^nC_{ij} + \sum_{i=0}^M\sum_{i=0}^M u_i^n u_j^n A_{ij}\\
 &= \mathbf{\hat{u}}^T \hat{A} \mathbf{\hat{u}} - 2 \mathbf{\hat{u}}^T C \mathbf{u}^n + \mathbf{u}^n A \mathbf{u}^n,
\end{align*}
where $A_{ij} = \langle \varphi_i,\varphi_j \rangle_{L^2}$, $\hat{A}_{ij} = \langle \hat{\varphi}_i,\hat{\varphi}_j \rangle_{L^2}$ and $C_{ij} = \langle \hat{\varphi}_i,\varphi_j \rangle_{L^2}$. 
The problem (\ref{eq:minprob}) can thus be reformulated as 
\begin{align*}
\min_{\hat{\mathbf{u}} \in \mathbb{R}^{\hat{M}+1}} \mathbf{\hat{u}}^T \hat{A} \mathbf{\hat{u}} - 2 \mathbf{\hat{u}}^T C \mathbf{u}^n + \mathbf{u}^n A \mathbf{u}^n \quad \text{s.t.} \quad \mathcal{I}_{\mathbf{p}^{n+1}}(\mathbf{\hat{u}}) - \mathcal{I}_{\mathbf{p}^{n}}(\mathbf{u}^n) = 0.
\end{align*}
This is a quadratic minimization problem with one nonlinear equality constraint, for which the solution $\hat{\mathbf{u}}$ is the solution of the nonlinear system of equations
\begin{align*}
\hat{A}\hat{\mathbf{u}} - C \mathbf{u}^n - \lambda \nabla \mathcal{I}_{\mathbf{p}^{n+1}}(\hat{\mathbf{u}}) &= 0\\
\mathcal{I}_{\mathbf{p}^{n+1}}(\mathbf{\hat{u}}) - \mathcal{I}_{\mathbf{p}^{n}}(\mathbf{u}^n) &= 0,
\end{align*}
which can be solved numerically using a suitable nonlinear solver.

\section{ADAPTIVE ENERGY PRESERVING SCHEMES FOR THE BBM EQUATION}

\subsection{The BBM equation}

The BBM equation was introduced by Peregrine \cite{peregrine66}, and later studied by Benjamin et al. \cite{BBM72} as a model for small amplitude long waves on the surface of water in a channel.
Conservative finite difference schemes for the BBM equation were proposed in \cite{wang06} and \cite{cohen11gfd}, the latter being a discrete gradient method on fixed grids. A moving mesh FEM scheme employing a quasi-Lagrange approach is presented by Lu, Huang and Qiu in \cite{Lu_Huang_Qiu}, which we also refer to for a more extensive list of references to the existing numerical schemes for the BBM equation.

Consider now an initial-boundary value problem of the one-dimensional BBM equation with periodic boundary conditions,
\begin{align}
&u_t - u_{x x t} + u_x + u u_{x} = 0,  &&\quad x \in [-L,L], \quad t \in (0,T] \label{eq:bbm}\\
&u(x,0) = u_0(x), &&\quad x \in [-L,L] \\
&u(-L,t) = u(L,t), &&\quad t \in (0,T]. \label{eq:periodic}
\end{align}
By introducing the new variable $m(x,t) := u(x,t) - u_{x x}(x,t)$, equation (\ref{eq:bbm}) can be rewritten on the form (\ref{eq:pde_var}) as
\begin{equation*}
m_t = \mathcal{S}(m) \frac{\delta \mathcal{H}}{\delta m},
\end{equation*}
for two different pairs of an antisymmetric differential operator $\mathcal{S}(m)$ and a Hamiltonian $\mathcal{H}\left[m\right]$:
\begin{eqnarray*}
\mathcal{S}^1(m) &=& -(\frac{2}{3}u+1)\partial_x - \frac{1}{3}u_x, \nonumber\\
\mathcal{H}^1\left[m\right] &=& \frac{1}{2} \int (u^2+u_x^2)\mathrm{d}x,
\end{eqnarray*}
and
\begin{eqnarray*}
\mathcal{S}^2(m) &=& -\partial_x + \partial_{x x x}, \nonumber\\
\mathcal{H}^2\left[m\right] &=& \frac{1}{2} \int (u^2+\frac{1}{3}u^3)\mathrm{d}x.
\end{eqnarray*}

\subsection{Discrete schemes}
We apply the PUM approach to create numerical schemes which preserve an approximation to either $\mathcal{H}^1\left[m\right]$ or $\mathcal{H}^2\left[m\right]$, splitting $\Omega := [-L,L]$ into $M$ elements $\{[x_i, x_{i+1}]\}_{i=0}^{M-1}$. Defining the matrices $A$ and $E$ by their components
\begin{align*}
A_{ij} =  \int_{\Omega} \varphi_{i}\varphi_{j} \mathrm{d}x \quad \mathrm{ and } \quad E_{ij} =  \int_{\Omega} \varphi_{i,x}\varphi_{j,x} \mathrm{d}x,
\end{align*}
we set $\mathbf{m} = (A+E)\mathbf{u}.$ Note that the matrices $A$ and $E$ depend on the mesh, and thus will change when adaptivity is used. We will then distinguish between matrices from different time steps by writing e.g. $A^n$ and $A^{n+1}$.

Approximating $u$ by $u^h$ as in section \ref{sec:PUM}, we find 
\begin{align*}
\mathcal{H}^1_{\mathbf{p}}(\mathbf{m}) &= \mathcal{H}^1[m^h] = \frac{1}{2} \int_{\Omega} (u^h)^2 + (u^h_x)^2 \mathrm{d}x \nonumber \\
&= \frac{1}{2} \sum \limits_{i,j} u_i u_j \int_{\Omega} \varphi_{i}\varphi_{j} \mathrm{d}x + \frac{1}{2} \sum \limits_{i,j} u_i u_j \int_{\Omega} \varphi_{i,x}\varphi_{j,x} \mathrm{d}x \\
&= \frac{1}{2}\mathbf{u}^\text{T}(A+E)\mathbf{u} \nonumber
\end{align*}
The integrals can be evaluated exactly and efficiently by considering elementwise which basis functions are supported on the element before applying Gaussian quadrature to obtain exact evaluations of the polynomial integrals. We define the matrix $B_1(\mathbf{u})$ by
\begin{align*}
B_1(\mathbf{u})_{j i} = -\frac{2}{3}\sum_{k=0}^{M-1} u_k \int_{\Omega} \varphi_{i}\varphi_{j,x} \varphi_{k} \mathrm{d}x - \int_{\Omega} \varphi_{i}\varphi_{j,x} \mathrm{d}x - \frac{1}{3}\sum_{k=0}^{M-1} u_k \int_{\Omega} \varphi_{i}\varphi_{j} \varphi_{k,x} \mathrm{d}x.
\end{align*} 
An approximation to the gradient of $\mathcal{H}^1$ with respect to $m$ is found by the AVF discrete gradient
\begin{align*}
\overline{\nabla}\mathcal{H}^1_{\mathbf{p}}(\mathbf{m}^n,\mathbf{m}^{n+1}) &= (A+E)^{-1}\overline{\nabla}\mathcal{H}^1_{\mathbf{p}}(\mathbf{u}^n,\mathbf{u}^{n+1}) \\
&= (A+E)^{-1}\int \limits_{0}^{1} \nabla\mathcal{H}^1_{\mathbf{p}}(\xi \mathbf{u}^{n} + (1-\xi)\mathbf{u}^{n+1}) \mathrm{d}\xi\\
&=(A+E)^{-1}\frac{1}{2}(A+E)\left(\mathbf{u}^{n}+\mathbf{u}^{n+1}\right) = \frac{1}{2}\left(\mathbf{u}^{n}+\mathbf{u}^{n+1}\right).
\end{align*}

Thus we have the required terms for forming the system (\ref{eq:PUM_ODE}) and applying the adaptive discrete gradient method to it. Corresponding to (\ref{eq:DGMM}), we get the scheme
\begin{align*}
(A^{n+1}+E^{n+1})\left(\mathbf{u}^{n+1} - \mathbf{\hat{u}}\right) = &\frac{\left(\mathbf{\hat{u}}^\text{T}\left(A^{n+1}+E^{n+1}\right)\mathbf{\hat{u}} - \left(\mathbf{u}^n\right)^\text{T}\left(A^{n}+E^{n}\right)\mathbf{u}^n\right)\left(\mathbf{\hat{u}}+\mathbf{u}^{n+1}\right)}{\left(\mathbf{\hat{u}}+\mathbf{u}^{n+1}\right)^\text{T}\left(\mathbf{\hat{u}}+\mathbf{u}^{n+1}\right)}\\
&+ \frac{\Delta t}{2}B^{n+1}_1\left(\frac{\mathbf{\hat{u}}+\mathbf{u}^{n+1}}{2}\right) \left(\mathbf{\hat{u}}+\mathbf{u}^{n+1}\right).
\end{align*}
Here we have chosen the skew-symmetric matrix $B_1$ to be a function of $\mathbf{\hat{u}}$ and $\mathbf{u}^{n+1}$, but could also have chosen e.g. $B_1(\mathbf{\hat{u}})$, resulting in a decreased computational cost at the expense of less precise results. During testing, the basis functions were chosen as piecewise cubic polynomials.

In the same manner we may obtain a scheme that preserves $\mathcal{H}^2\left[m\right]$. In this case
\begin{align*}
\mathcal{H}^2_{\mathbf{p}}(\mathbf{m}) &= \mathcal{H}^2[m^h]  = \frac{1}{2} \int_{\Omega} (u^h)^2  + \frac{1}{3}(u^h)^3 \mathrm{d}x \nonumber \\
&=  \frac{1}{2} \sum \limits_{i,j} u_i u_j \int_{\Omega} \varphi_{i}\varphi_{j} \mathrm{d}x + \frac{1}{6} \sum \limits_{i,j,k}   u_i u_j u_k \int_{\Omega} \varphi_{i}\varphi_{j}\varphi_{k} \mathrm{d}x.
\end{align*}
and
\begin{align*}
(B_2)_{ji} = - \int_{\Omega} \varphi_{i}\varphi_{j,x} \mathrm{d}x + \int_{\Omega} \varphi_{i}\varphi_{j,xxx} \mathrm{d}x.
\end{align*}
Note that the skew-symmetric matrix $B_2$ is independent of $\mathbf{u}$.

Defining the tensor $D$ by its elements
\begin{align*}
D_{ijk} = \int_{\Omega} \varphi_{i}\varphi_{j}\varphi_{k} \mathrm{d}x,
\end{align*}
we get, with the convention of summation over repeated indices, the AVF discrete gradient with respect to $\mathbf{u}$ given by the elements
\begin{align*}
\left(\overline{\nabla}\mathcal{H}^2_{\mathbf{p}}(\mathbf{u}^n,\mathbf{u}^{n+1})\right)_i  = \dfrac{1}{2} A_{ij}(u_j^n + u_j^{n+1}) + \frac{1}{6} D_{ijk}(u_j^n (u_k^n  + \frac{1}{2} u_k^{n+1}) + u_j^{n+1}(\frac{1}{2} u_k^n  + u_k^{n+1})).
\end{align*}
and again the discrete gradient with respect to $\mathbf{m}$ by
$$\overline{\nabla}\mathcal{H}^2_{\mathbf{p}}(\mathbf{m}^n,\mathbf{m}^{n+1}) = (A+E)^{-1}\overline{\nabla}\mathcal{H}^2_{\mathbf{p}}(\mathbf{u}^n,\mathbf{u}^{n+1}).$$
If we employ integral preserving transfer between the meshes, we get the scheme
\begin{align*}
\mathbf{u}^{n+1} - \mathbf{\hat{u}} = \Delta t (A+E)^{-1}B_2 (A+E)^{-1} \overline{\nabla}\mathcal{H}^2_{\mathbf{p}}(\mathbf{\hat{u}},\mathbf{u}^{n+1}),
\end{align*}
where we note that $S_{\mathbf{p},2}:=(A+E)^{-1}B_2 (A+E)^{-1}$ is a skew-symmetric matrix. If non-preserving transfer is used, we need a correction term, as in the $\mathcal{H}^1$ scheme above. The calculation of such a term is straightforward, but we omit it here for reasons of brevity.

To approximate the third derivative in $B_2$, we need basis functions of at least degree three, and to guarantee skew-symmetry in $B_2$, these basis functions need to be $C^2$ on the element boundaries. This is not obtainable with regular nodal FEM basis functions, so we have instead used third order B-spline basis functions as described in \cite{Hughes} during testing.

\section{NUMERICAL RESULTS}
To demonstrate the performance of our methods, we have tested them on two one-dimensional simple problems: A soliton solution, and the interaction of two waves. We have tested our $\mathcal{H}^1$- and $\mathcal{H}^2$-preserving schemes on uniform and moving meshes, and compared the results to those obtained using the explicit midpoint method. For the transfer operation between meshes, we have used a piecewise cubic interpolation method in the $\mathcal{H}^1$ preserving scheme, and exact transfer in the $\mathcal{H}^2$ preserving scheme.

\subsection{Mesh adaptivity}
As noted in section \ref{sec:adaptive}, our methods can be coupled with any adaptive mesh strategy using the rezoning approach. For our numerical experiments, we have used a simple method for $r$-adaptivity based on the equidistribution principle: Splitting $\Omega$ into $M$ intervals, we require that
\begin{align*}
\int \limits_{x_i}^{x_{i+1}} \omega(x) \mathrm{d}x = \frac{1}{M} \int \limits_{-L}^{L} \omega(x) \mathrm{d}x,
\end{align*}
where the monitor function $\omega$ is a function measuring how densely grid points should lie, based on the value of $u$. For a general discussion on the choice of an optimal monitor function, see e.g. \cite{BuddHuangRussell, BlomVerwer}. For the problems we have studied, a generalized solution arc length monitor function proved to yield good results. This is given by
\begin{align*}
\omega(x) = \sqrt{1 + k^2\left( \dfrac{\partial u}{\partial x}(x) \right)^2}.
\end{align*}
For $k=1$, this is the usual arc length monitor function, in which case the equidistribution principle amounts to requiring that the arc length of $u$ over each interval is equal. In applications, we only have an approximation of $u$, and hence $\omega$ must be approximated as well. We have applied a finite difference approximation and obtained approximately equidistributing grids using de Boor's method as explained in \cite[pp.~36-38]{MM_Huang_Russell}.

\subsection{Soliton solution}

With $u_0(x) = 3(c-1)\hspace{2pt}\text{sech}^2\left(\frac{1}{2}\sqrt{1-\frac{1}{c}}x\right)$, the exact solution of (\ref{eq:bbm})--(\ref{eq:periodic}) is
\begin{align*}
u(x,t) = 3(c-1)\hspace{1pt}\text{sech}^2\left(\frac{1}{2}\sqrt{1-\frac{1}{c}} l(x,t)\right),
\end{align*}
with $l(x,t) = \min_{j \in \mathbb{Z}} \left|x-c t+2 j L\right|$.
This is a soliton solution which travels with a constant speed $c$ in $x$-direction while maintaining its initial shape.

To evaluate the numerical solutions, we have compared them to the exact solution and calculated errors in shape and phase. The phase error is evaluated as
\begin{align*}
E^{\text{phase}}_n = |ct_n - x^*|,
\end{align*}
where $x^* = \text{arg} \max \limits_{x} u_h(x,t_n) $, i.e. the location of the peak of the soliton in the numerical solution. The shape error is given by
\begin{align*}
E^{\text{shape}}_n = \left|\left|u_h(x,t_n) - u \left(x,\frac{x^*}{c} \right) \right| \right|,
\end{align*}
where the peak of the exact solution is translated to match the peak of the numerical solution, and the difference in the shapes of the solitons is calculated.

\begin{figure}
        \centering
        \subfloat{
                \centering
                \includegraphics[width=0.49\textwidth]{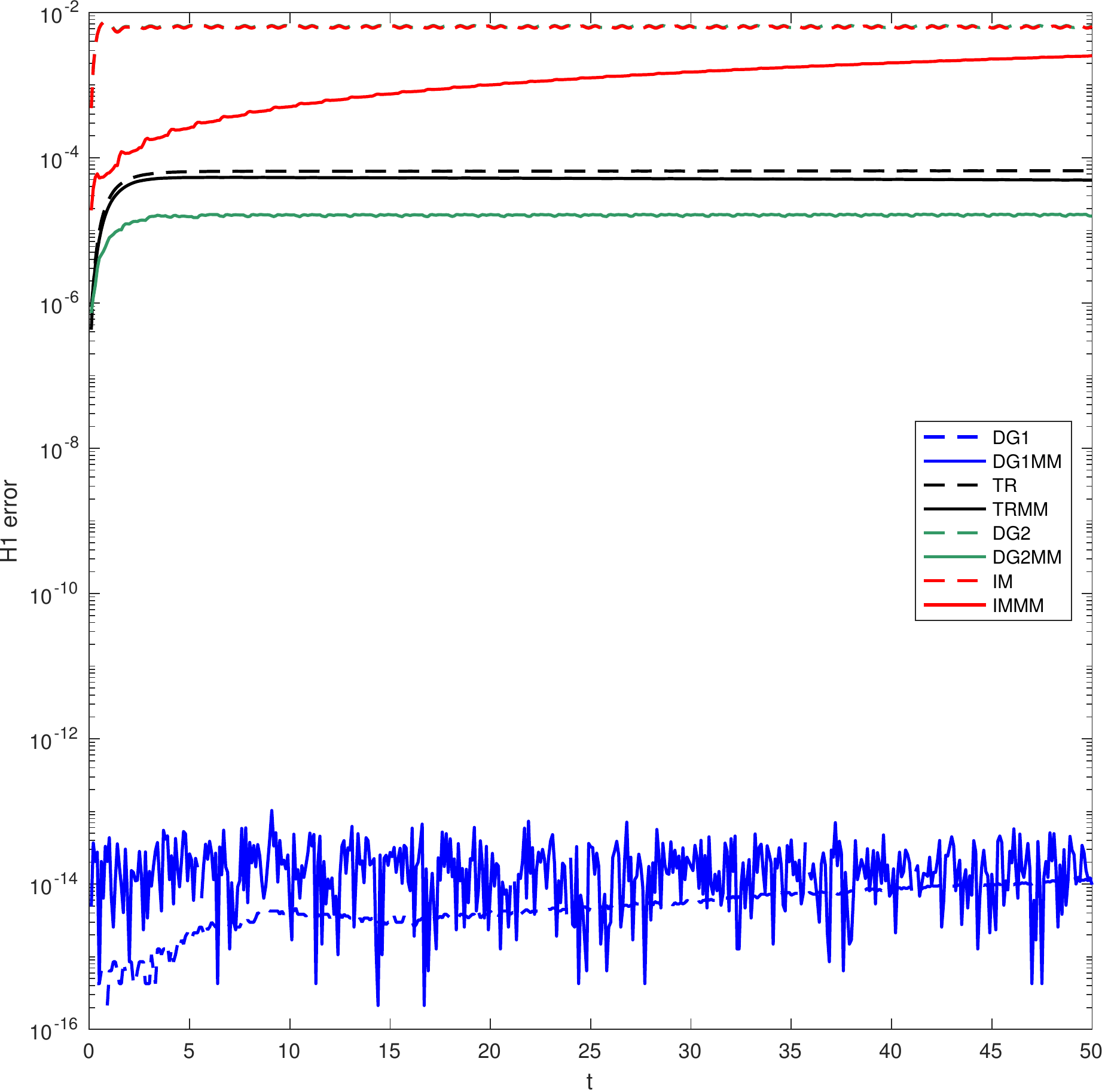}
        }
        \subfloat{
                \centering
                \includegraphics[width=0.49\textwidth]{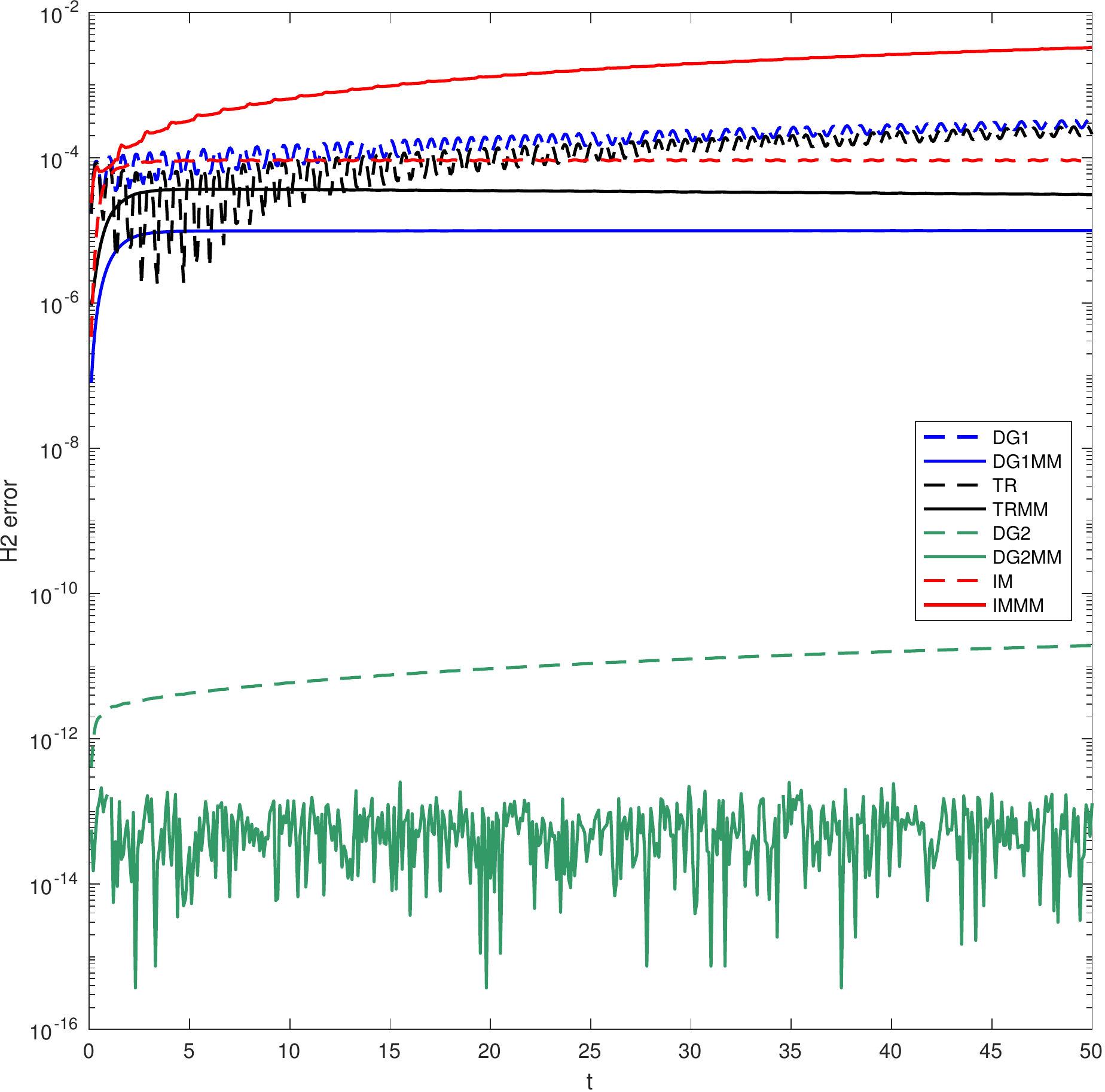}
        }
        \caption{The soliton problem. Relative error in the approximated Hamiltonians $\mathcal{H}^1_{\mathbf{p}}$ (left) and $\mathcal{H}^2_{\mathbf{p}}$ (right) plotted as a function of time $t \in \left[0,50\right]$. $c = 3, L = 200, \Delta t = 0.1$, $M = 200$.}
        \label{fig:energyerrors}
\end{figure}

\begin{figure}
        \centering
        \subfloat{
                \centering
                \includegraphics[width=0.49\textwidth]{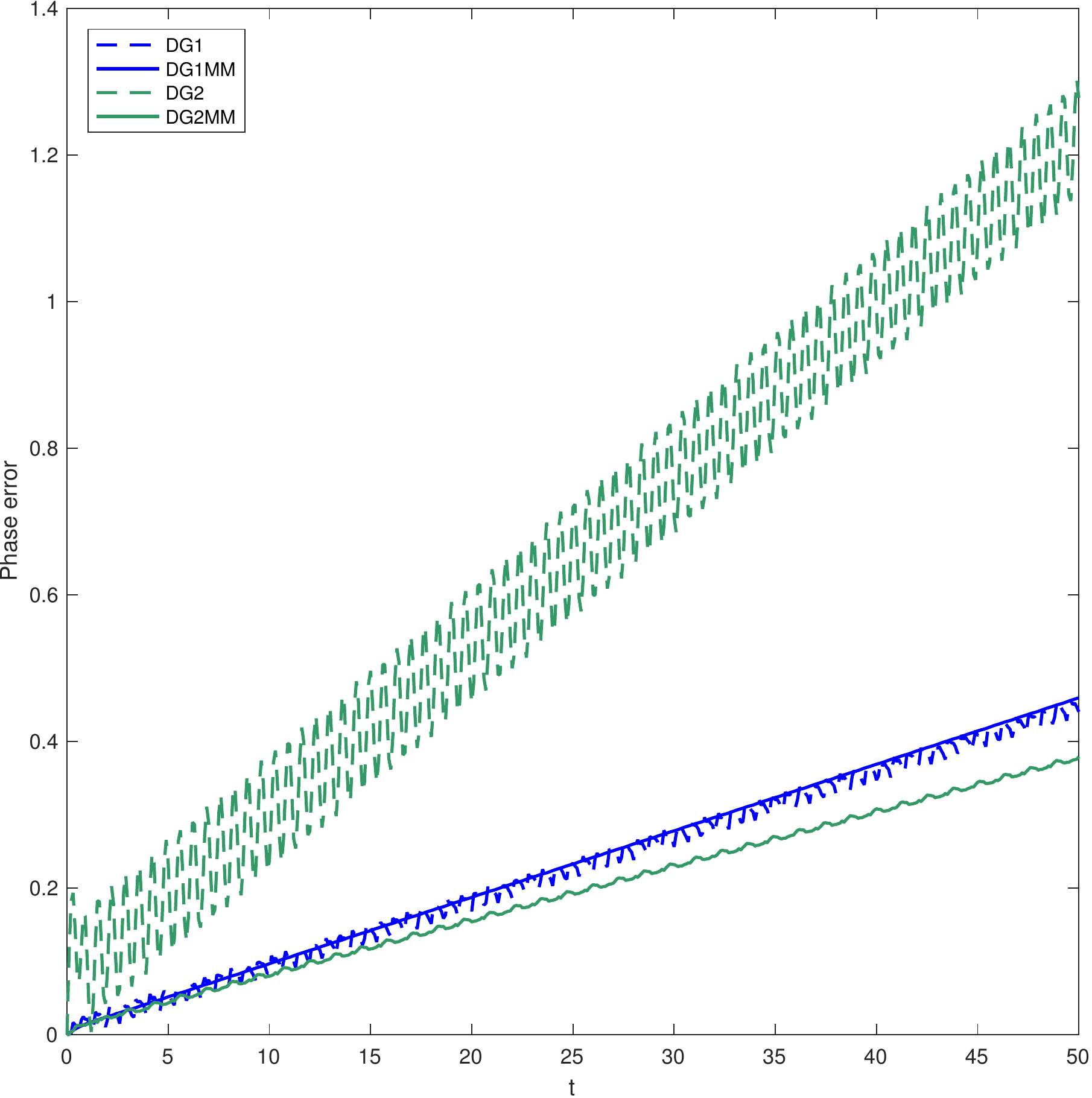}
        }
        \subfloat{
                \centering
                \includegraphics[width=0.49\textwidth]{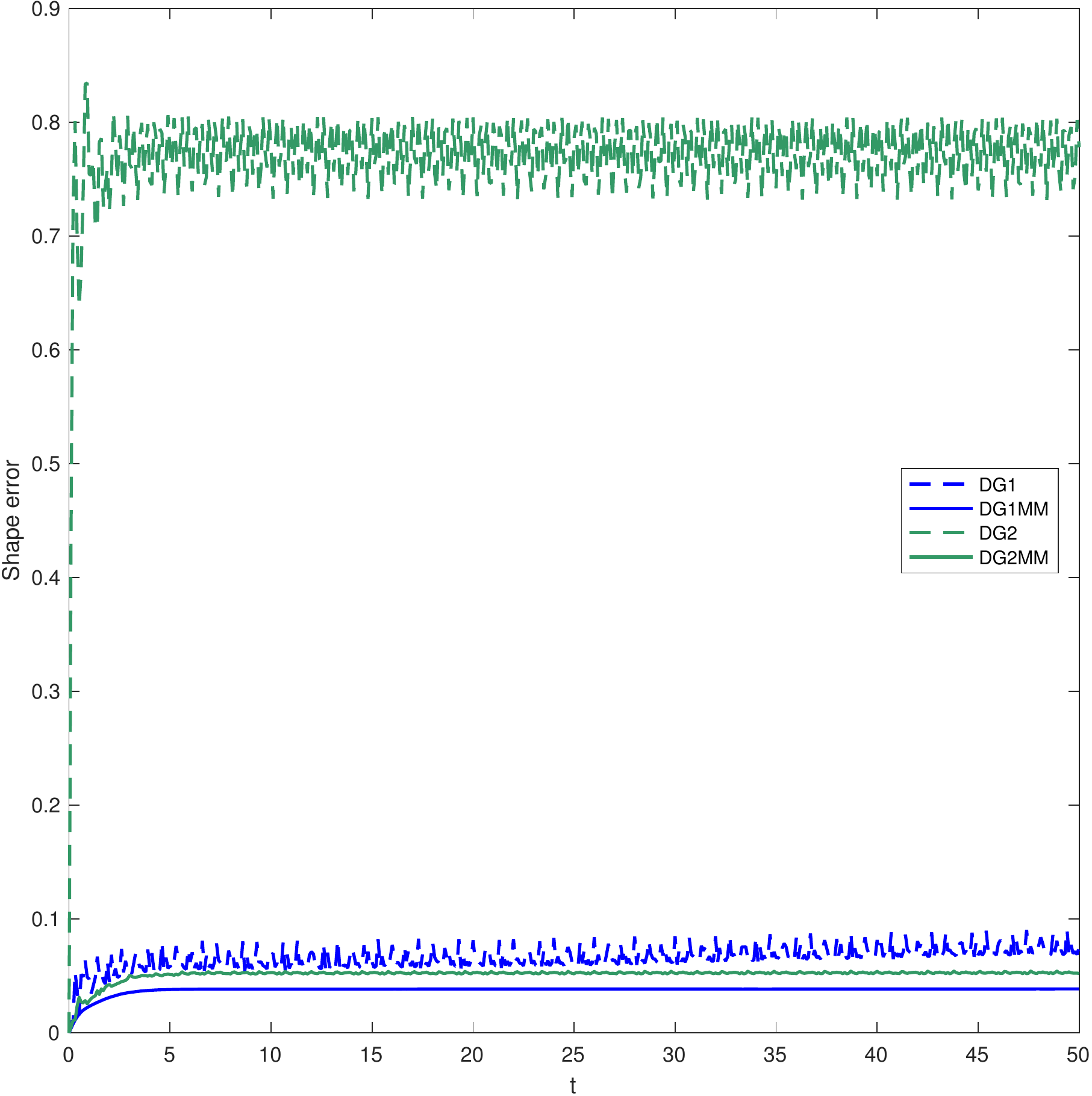}
        }
        \caption{The soliton problem. Phase error (left) and shape error (right) as a function of time. $c = 3, L = 200, \Delta t = 0.1$, $M = 200$.}
        \label{fig:overtime_mm}
\end{figure}

\begin{figure}
        \centering
        \subfloat{
                \centering
                \includegraphics[width=0.49\textwidth]{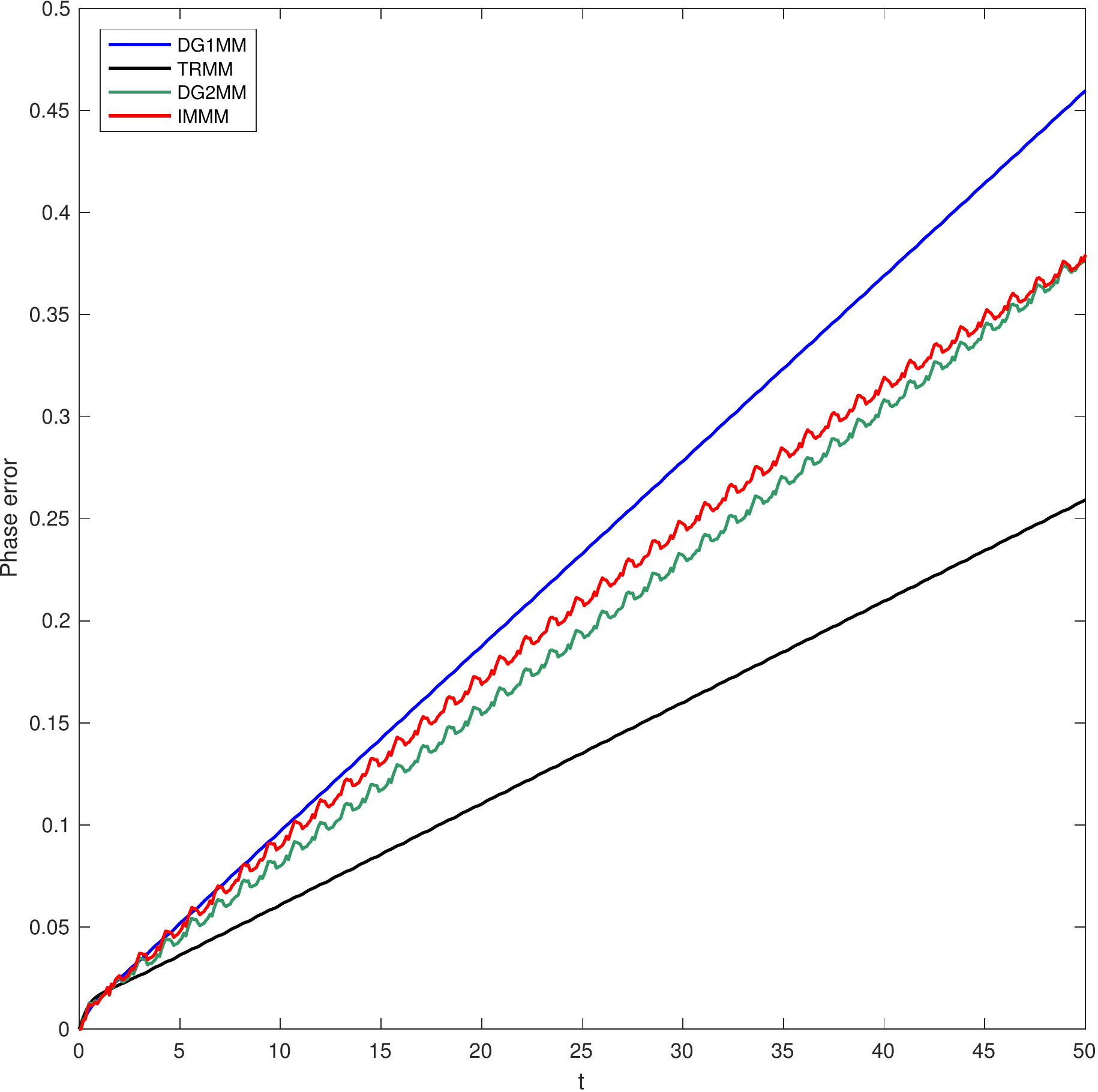}
        }
        \subfloat{
                \centering
                \includegraphics[width=0.49\textwidth]{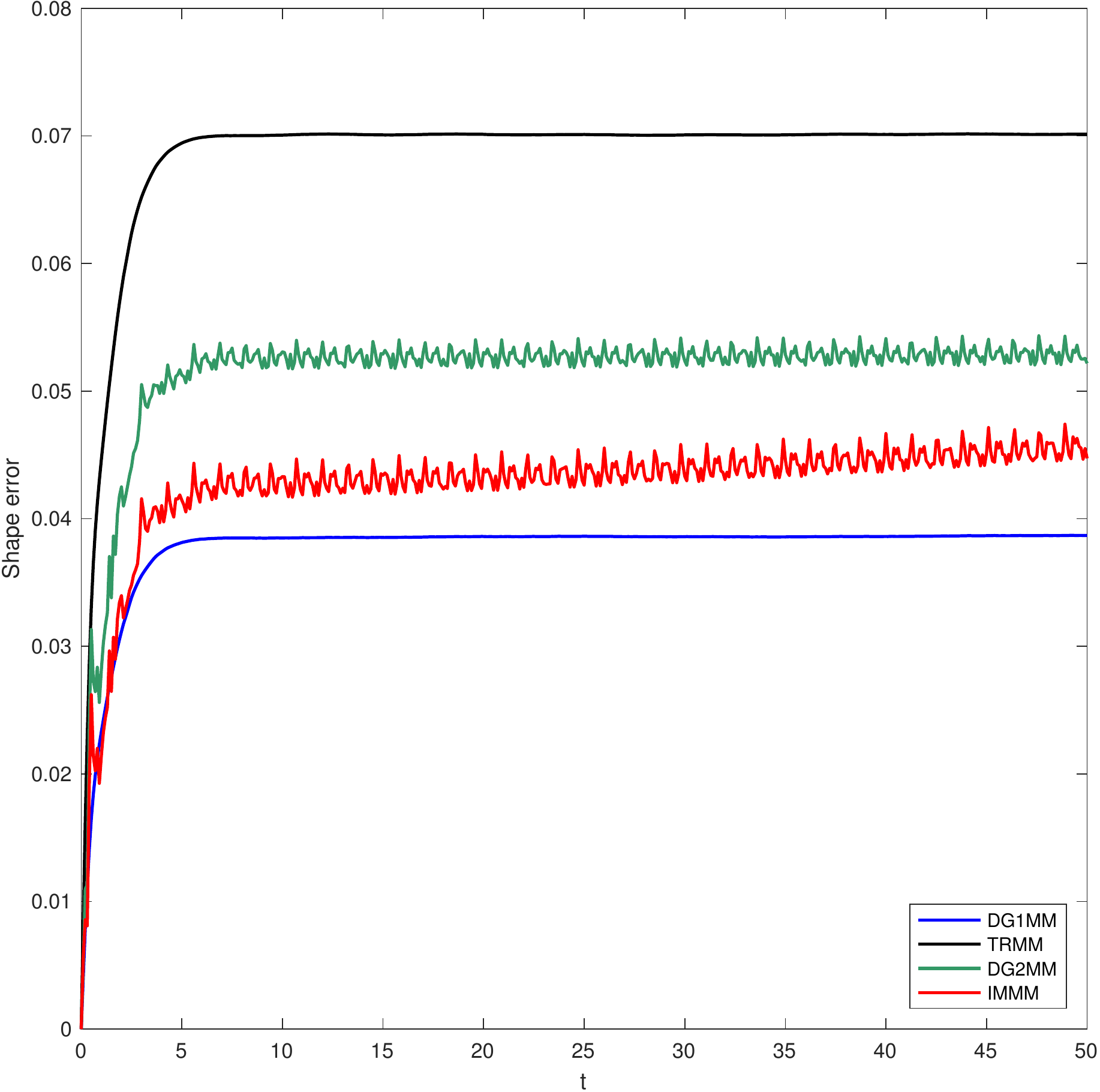}
        }
        \caption{The soliton problem. Phase error (left) and shape error (right) as a function of time. $c = 3, L = 200, \Delta t = 0.1$, $M = 200$.}
        \label{fig:overtime_dg}
\end{figure}

The results of the numerical tests can be seen in figures \ref{fig:energyerrors}--\ref{fig:overtime_dg}. Here, $M$ denotes the degrees of freedom used in the spatial approximation and $\Delta t$ the fixed time step size. DG1 and DG1MM denotes the $\mathcal{H}^1_{\mathbf{p}}$ preserving scheme with fixed, uniform grid and adaptive grid, respectively; similary DG2 and DG2MM denotes the $\mathcal{H}^2_{\mathbf{p}}$ preserving scheme with uniform and adaptive grids.

In Figure \ref{fig:energyerrors} we see the relative errors in $\mathcal{H}^1_{\mathbf{p}}$ and $\mathcal{H}^2_{\mathbf{p}}$. The DG1 and DG1MM schemes are compared to schemes using the same 3rd order nodal basis functions, but the trapezoidal rule for time-stepping, denoted by TR and TRMM. Likewise, the DG2 and DG2MM schemes are compared to the IM and IMMM schemes, using B-spline basis functions and the implicit midpoint method for discretization in time. The error in $\mathcal{H}^1_{\mathbf{p}}$ is very small for the DG1 and DG1MM schemes, as expected. Also the error in $\mathcal{H}^2_{\mathbf{p}}$ is very small for the DG2 and DG2MM schemes. The order of the error is not machine precision, but is instead dictated by the precision with which the nonlinear equations in each time step is solved. We can also see that while the TR and IM schemes, with and without moving meshes, have poor conservation properties, the moving mesh DG schemes seem to preserve quite well even the integrals they are not designed to preserve.

In figures \ref{fig:overtime_mm} and \ref{fig:overtime_dg} we see the phase and shape errors, of our methods compared to non-moving mesh methods and non-preserving methods, respectively. The advantage of using moving meshes is clear, especially for the $\mathcal{H}^2_{\mathbf{p}}$ preserving schemes. The usefulness on integral preservation is ambiguous in this case. It seems that what we gain in precision in phase, we lose in precision in shape, and vice versa.

\subsection{A small wave overtaken by a large one}
A typical test problem for the BBM equation is the interaction between two solitary waves. With an initial condition
$$u_0(x) = 3(c_r-1)\hspace{2pt}\text{sech}^2\left(\frac{1}{2}\sqrt{1-\frac{1}{c_r}}(x-x_r)\right) + 3(c_s-1)\hspace{2pt}\text{sech}^2\left(\frac{1}{2}\sqrt{1-\frac{1}{c_s}}(x-x_s)\right),$$
one wave will eventually be overtaken by the other as long as $c_r \neq c_s$, i.e. if one wave is larger than the other. There is no available analytical solution for this problem. The two waves are not solitons, as the amplitudes will change a bit after the waves have interacted \cite{craig2006}.

\begin{figure}
        \centering
        \subfloat{
                \centering
                \includegraphics[width=0.49\textwidth]{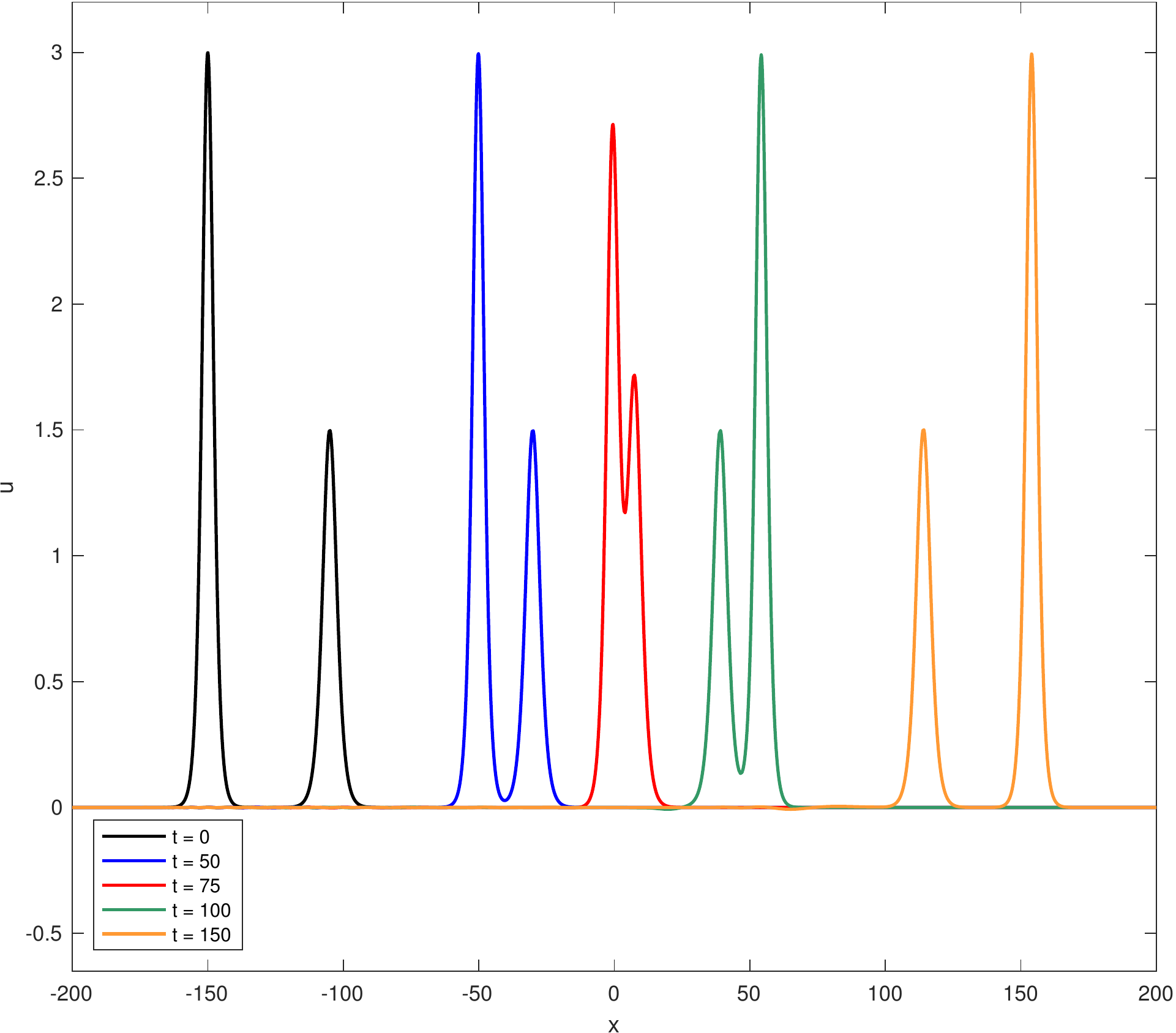}
        }
        \subfloat{
                \centering
                \includegraphics[width=0.49\textwidth]{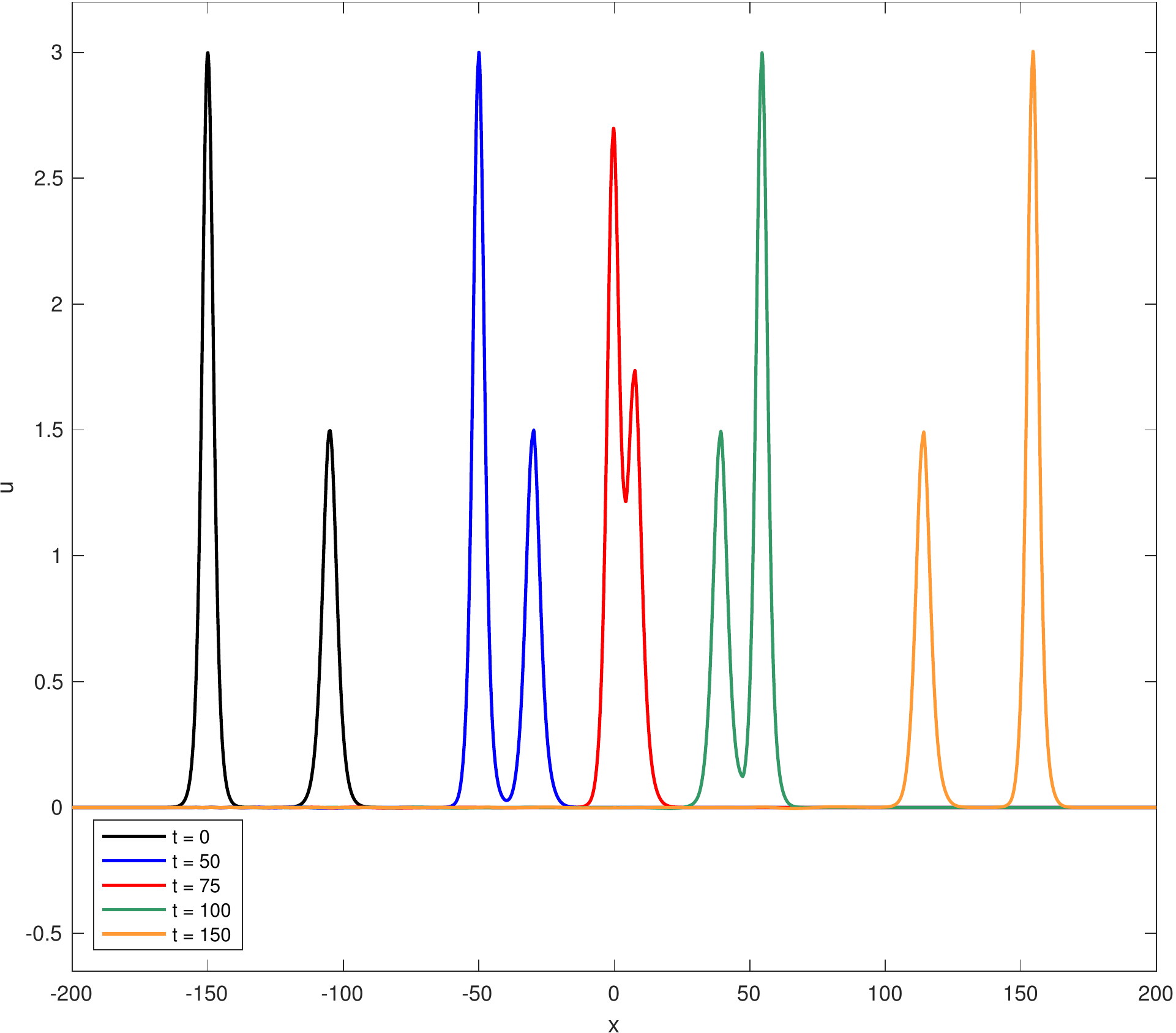}
        }
        \caption{The interacting waves problem. Solutions at $t = \left\{0,50,75,100,150\right\}$ found by DG1MM (left) and DG2MM (right). $x_r = 150, x_s = 105, c_r = 2, c_s = 1.5, L = 200, \Delta t = 0.1$, $M = 1000$.}
        \label{fig:2waves}
\end{figure}

\begin{figure}
        \centering
                \includegraphics[width=0.7\textwidth]{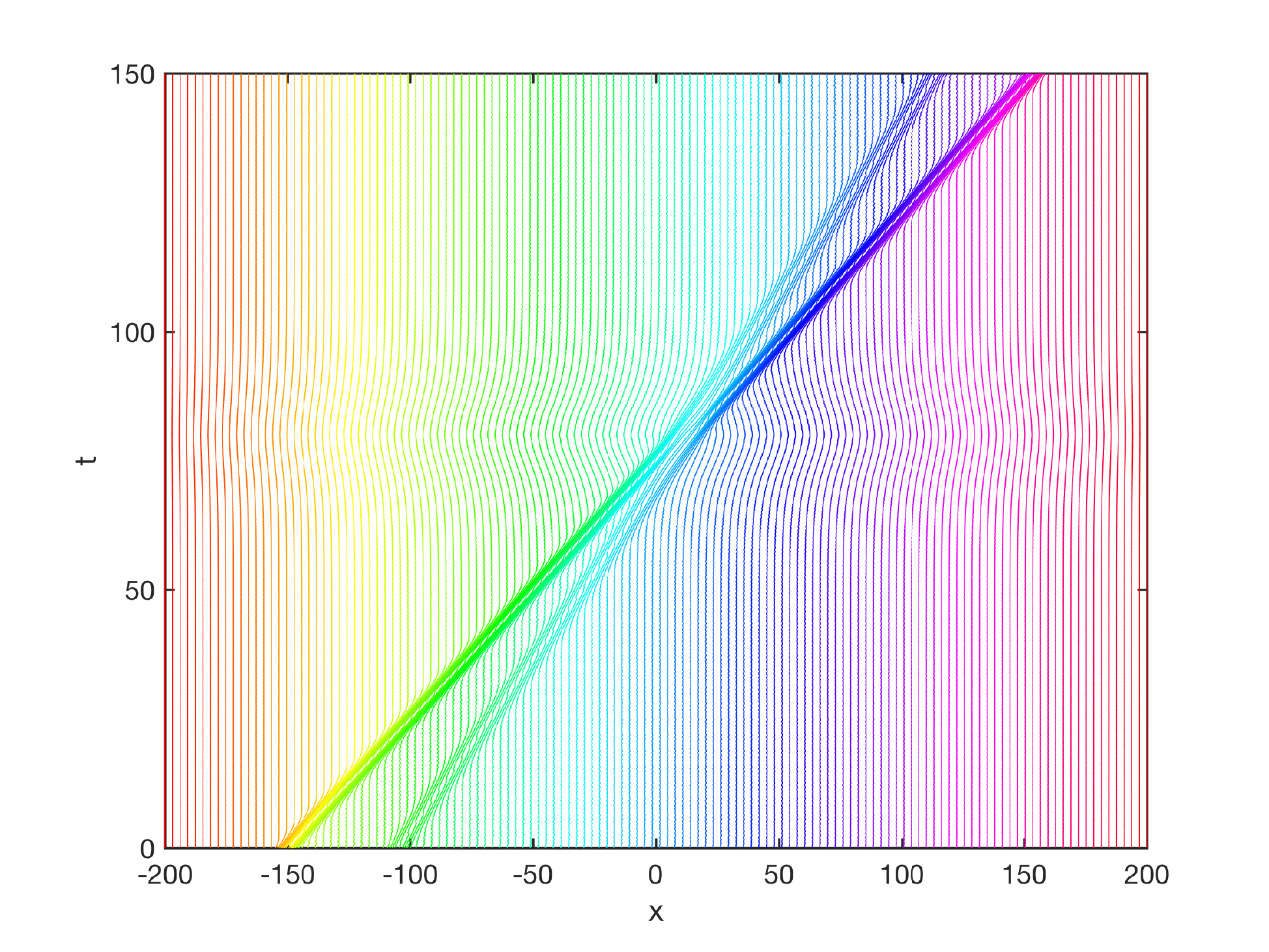}
        \caption{Mesh point trajectories in time. Each line represents one mesh point.}
        \label{fig:mmplot}
\end{figure}

Solutions obtained by solving the problem with our two energy preserving schemes, giving very similar results, are plotted in Figure \ref{fig:2waves}. Also, to illustrate the mesh adaptivity, we have included a plot of the mesh trajectories in Figure \ref{fig:mmplot}. Each line represents the trajectory of one mesh point in time, and we can see that the mesh points cluster nicely around the edges of the waves as they move. 

\begin{figure}
        \centering
        \subfloat{
                \centering
                \includegraphics[width=0.49\textwidth]{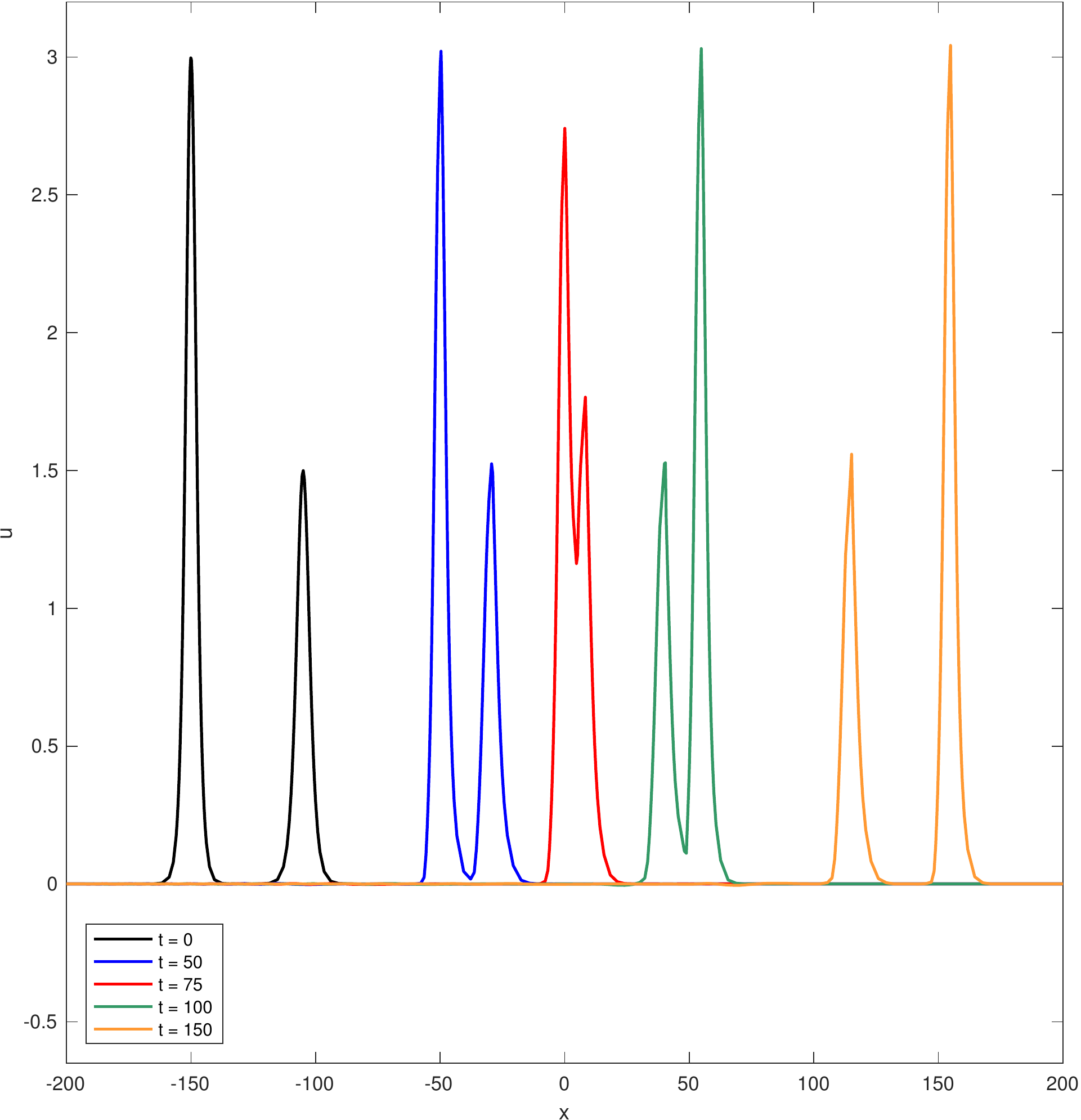}
        }
        \subfloat{
                \centering
                \includegraphics[width=0.49\textwidth]{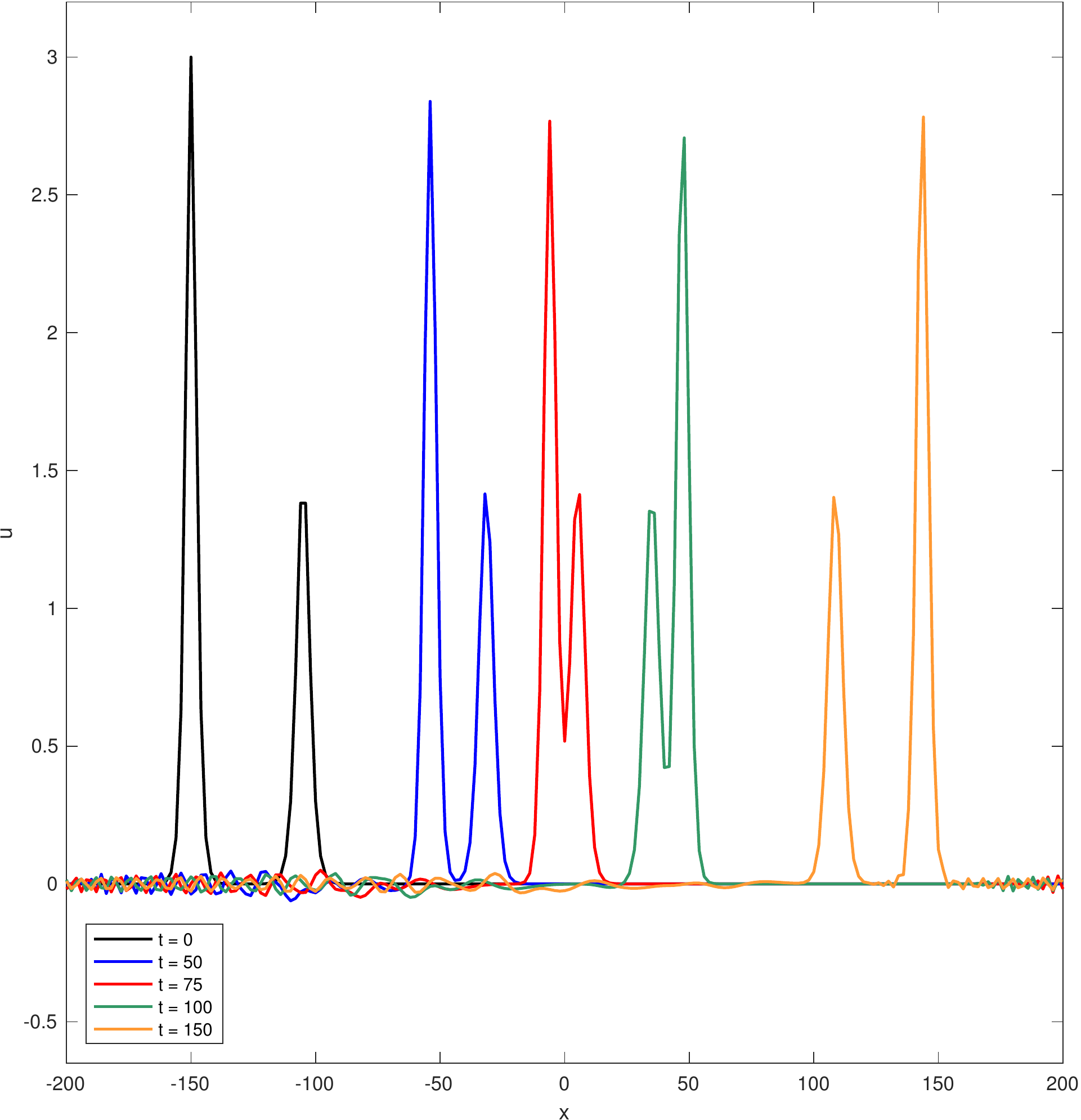}
        }
        \caption{The interacting waves problem. Solutions at $t = \left\{0,50,75,100,150\right\}$ found by DG2MM (left) and RK (right). $x_r = 150, x_s = 105, c_r = 2, c_s = 1.5, L = 200, \Delta t = 0.1$, $M = 200$.}
        \label{fig:2waves200}
\end{figure}

To illustrate the performance of our methods, we have in Figure \ref{fig:2waves200} compared solutions obtained by using the $\mathcal{H}_\mathbf{p}^2$-preserving moving mesh method with the solutions obtained by using a fourth order Runge--Kutta method on a static mesh, with the same, and quite few, degrees of freedom. The DG2MM solution is visibly closer to the solutions in Figure \ref{fig:2waves}. The non-preserving RK scheme does a worse job of preserving the amplitude and speed of the waves compared to the DG2MM scheme, and we observe unwanted oscillations. 

\begin{figure}
        \centering
        \subfloat{
                \centering
                \includegraphics[width=0.49\textwidth]{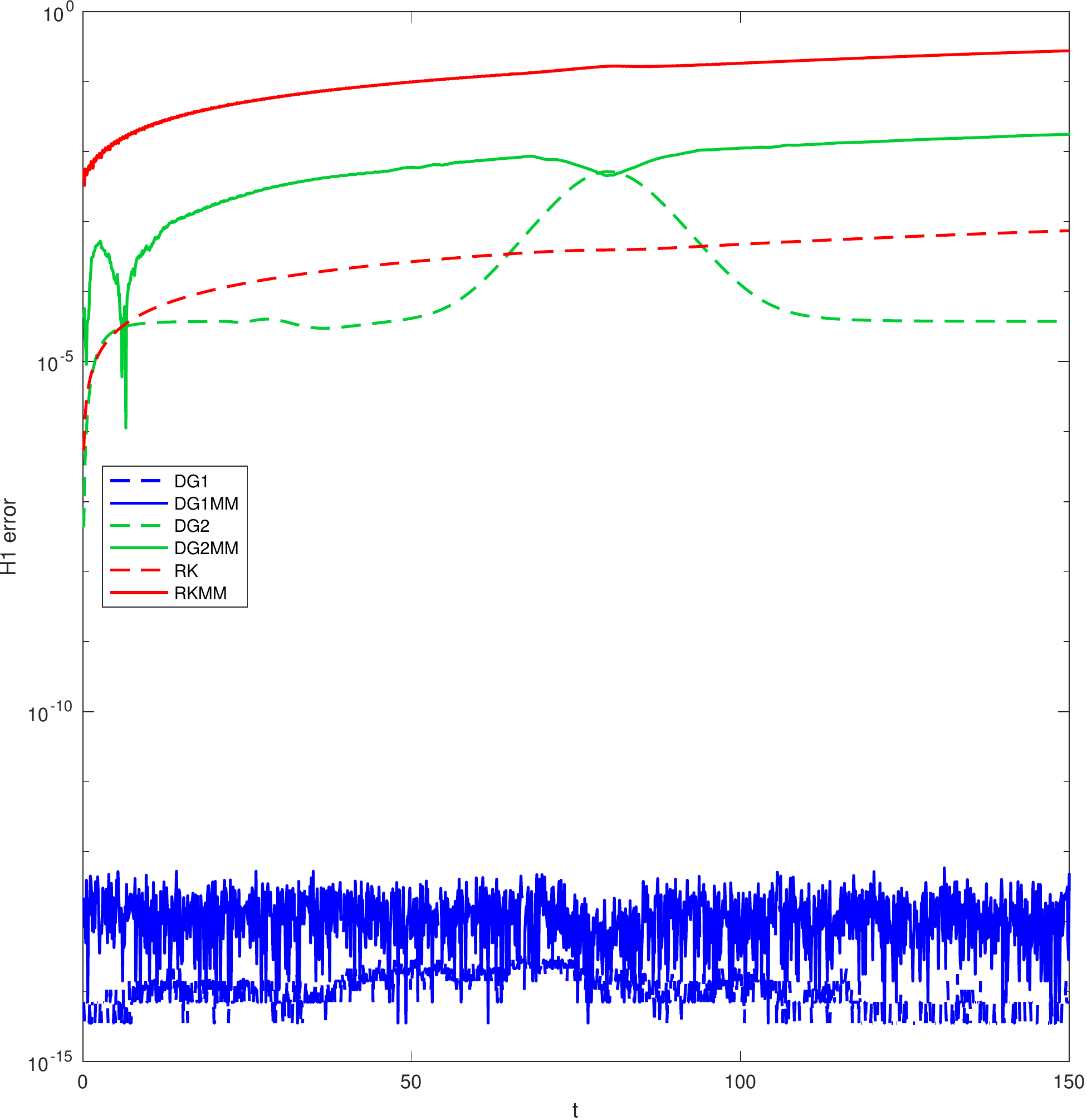}
        }
        \subfloat{
                \centering
                \includegraphics[width=0.49\textwidth]{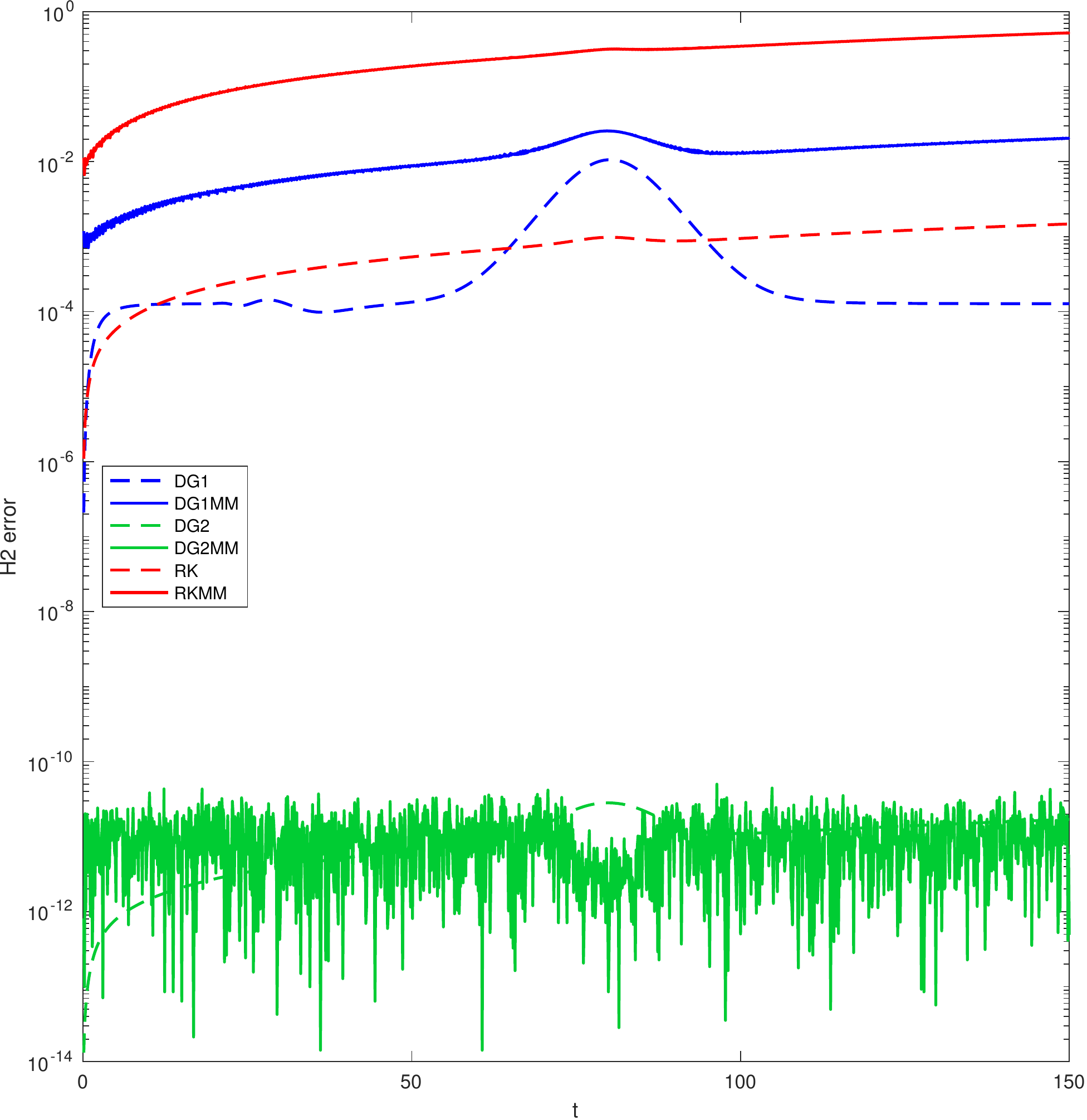}
        }
        \caption{The interacting waves problem. Error in the approximated Hamiltonians $\mathcal{H}^1_{\mathbf{p}}$ (left) and $\mathcal{H}^2_{\mathbf{p}}$ (right) plotted as a function of time $t \in \left[0,150\right]$. $x_r = 150, x_s = 105, c_r = 2, c_s = 1.5, L = 200, \Delta t = 0.1$, $M = 1000$.}
        \label{fig:energyerrors2w}
\end{figure}

In Figure \ref{fig:energyerrors2w} we have plotted the Hamiltonian errors for this problem. Again we see that the energy preserving schemes preserve both Hamiltonians better than the Runge--Kutta scheme, but we do also observe that the DG1 scheme preserves $\mathcal{H}^2_{\mathbf{p}}$ better than the DG1MM scheme, and vice versa for the DG2 and DG2MM schemes. Note also that an increase in the errors can be observed when the two waves interact, but that this increase is temporary.

\section{CONCLUSIONS}

In this paper, we have presented energy preserving schemes for a class of PDEs, first on general fixed meshes, and then on adaptive meshes. These schemes are then applied to the BBM equation, for which discrete schemes preserving two of the Hamiltonians of the problem are explicitly given.

Numerical experiments are performed, using the energy preserving moving mesh schemes on two different BBM problems: a soliton solution, and two waves interacting. Plots of the phase and shape errors illustrate how, for the given parameters, the usage of moving meshes gives improved accuracy, while the integral preservation gives comparable results to existing methods, without yielding a categorical improvement. We will remark, however, that in many cases, the preservation of a quantity such as one of the Hamiltonians in itself may be a desired property of a numerical scheme. For the two wave interaction problem, we do not have an analytical solution to compare to, but plots of the solution indicate that our schemes perform well compared to a Runge--Kutta scheme.

Although the numerical examples presented here are simple one-dimensional problems, the adaptive discrete gradient methods should also be applicable for multi-dimensional problems. This could be an interesting direction for further work, since the advantages of adaptive meshes are typically more evident when increasing the number of dimensions.

\bibliography{EP,tb}
\bibliographystyle{ieeetr}

\end{document}